\documentclass[11pt,a4paper]{article}
\usepackage[T1]{fontenc}
\usepackage[utf8]{inputenc}
\usepackage[margin=25mm]{geometry}
\usepackage{lmodern,microtype,amsmath,amssymb,amsthm,mathtools}
\usepackage{booktabs,longtable,array,enumitem,etoolbox,needspace,graphicx}
\usepackage[table]{xcolor}
\usepackage[colorlinks=true,linkcolor=blue!45!black,urlcolor=blue!45!black,citecolor=blue!45!black]{hyperref}
\usepackage{bookmark}

\usepackage{orcidlink}

\usepackage{xurl}
\usepackage{caption}
\usepackage{diagbox}
\usepackage{color}
\usepackage{placeins,float}
\definecolor{exactcell}{HTML}{E4F0E7}
\definecolor{boundcell}{HTML}{FAEADB}

\newcommand{\AuthorOneORCID}{0009-0005-6539-2249}
\newcommand{\AuthorTwoORCID}{0000-0002-6036-5541}
\newcommand{\AuthorOneLocation}{Reggio Calabria, Italy}
\newcommand{\AuthorTwoLocation}{Rome, Italy}
\newcommand{\AuthorOneEmail}{aldoroberto.pessolano@gmail.com}
\newcommand{\AuthorTwoEmail}{marcokrt1984@yahoo.it}
\newcommand{\AuthorORCID}[1]{\ifdefempty{#1}{}{\,\orcidlink{#1}}}
\newcommand{\AuthorLocation}[1]{\ifdefempty{#1}{}{\\\small #1}}
\newcommand{\AuthorEmail}[1]{\ifdefempty{#1}{}{\\\small\href{mailto:#1}{\expandafter\nolinkurl\expandafter{#1}}}}
\hypersetup{pdftitle={Covering hypercubic grids by polygonal trails: capacity bounds and the optimum for P(4,3)},pdfauthor={Aldo Roberto Pessolano; Marco Rip\`a}}

 \setlist{itemsep=3pt,topsep=5pt} \newtheorem{theorem}{Theorem}[section] \newtheorem{lemma}[theorem]{Lemma} \newtheorem{proposition}[theorem]{Proposition} \newtheorem{corollary}[theorem]{Corollary} \theoremstyle{definition}\newtheorem{definition}[theorem]{Definition}\newtheorem{problem}[theorem]{Problem} \theoremstyle{remark}\newtheorem{remark}[theorem]{Remark}  \newcommand{\card}[1]{\lvert #1\rvert}
\date{11 November 2026}
\title{Covering hypercubic grids by polygonal trails:\\capacity bounds and the optimum for $P(4,3)$}

\author{Aldo Roberto Pessolano\textsuperscript{1}\AuthorORCID{\AuthorOneORCID}\\
\small\textsuperscript{1}Independent Researcher\AuthorLocation{\AuthorOneLocation}\AuthorEmail{\AuthorOneEmail}
\and Marco Rip\`a\textsuperscript{2}\AuthorORCID{\AuthorTwoORCID}\\
\small\textsuperscript{2}Independent Researcher\AuthorLocation{\AuthorTwoLocation}\AuthorEmail{\AuthorTwoEmail}}
\begin{document}\maketitle
\begin{abstract}
The classical nine-dot puzzle asks for a polygonal trail through a square grid using as few segments as possible. We study its higher-dimensional analogue on $P(n,k):=\{0,\ldots,n-1\}^k$ for integers $n,k\ge2$, allowing arbitrary real vertices and unrestricted segment directions, and denote the minimum number of segments by $h(n,k)$. For $n\ge3$ and $k\ge2$, we give a geometric proof of the previously stated capacity bound: a trail of $r\ge1$ segments meets at most $1+r(n-1)$ grid points. Refinements of this bound, combined with exhaustive finite checks, establish $h(4,3)=23$ and $h(5,3)\ge32$. We also exhibit a $22$-segment trail through $63$ points of $P(4,3)$, the maximum possible, and a $23$-segment covering trail contained in $[0,4]^3$, which settles the corresponding box-constrained problem. The general capacity theorem, a large-grid lower bound refuting a previously proposed formula for $h(n,k)$, and the boxed construction are verified by the Lean 4 kernel, whereas the finite exclusions use native reflection and additionally trust the native compiler and runtime. A companion technical note provides the computational specification and reproducible proof artifacts.
\end{abstract}


\section{Introduction} \label{sec:Intro}
The classical nine-dot puzzle, appearing on p.~301 of Sam Loyd's posthumously published 1914 collection \cite{loyd}, asks for a continuous polygonal path through all the points of a $3\times3$ array using as few segments as possible. Its familiar four-segment solution uses turning points outside the array.

In 1970, Golomb and Selfridge \cite{golomb} gave a systematic treatment of minimum-segment polygonal paths and related closed constructions on planar point lattices. Later, Bereg et al. \cite{bereg} studied minimum-turn traversals with axis-aligned segments, while Keszegh \cite{keszegh} treated covering paths and trees for planar grids without that restriction.

The higher-dimensional analogue of the nine-dot puzzle replaces the square array by a hypercubic grid, with the number of points along each side and the dimension varying independently.

We use the following formulation, consistent with \cite{ripa}.
\begin{definition}[Grids and trails]\label{def:covering}
For integers $n\ge2$ and $k\ge2$, let
\[
 P(n,k):=\{0,1,\ldots,n-1\}^k\subset\mathbb R^k.
\]
A \emph{polygonal trail} of $r\ge1$ segments is specified by vertices $v_0,\ldots,v_r\in\mathbb R^k$ with $v_{i-1}\ne v_i$ and consecutive closed segments $[v_{i-1},v_i]$, $1\le i\le r$, such that segments in distinct positions have at most one point in common. This condition prohibits \emph{retracing}, namely traversing a subsegment of positive length more than once, in either direction. Isolated self-intersections and repeated visits to grid points or vertices are allowed. A vertex of the trail that does not belong to $P(n,k)$ is called a \emph{Steiner point}; it may lie inside or outside the smallest enclosing box $[0,n-1]^k$.
A trail is called \emph{closed} if its initial and final vertices coincide.
\end{definition}

The problem studied here is to minimize the number of segments in a covering trail. We denote the optimum by $h(n,k)$ and use $h_l(n,k)$ and $h_u(n,k)$ for a proved lower bound and a constructive upper bound, respectively.
\begin{problem}[Minimum-link covering]\label{prob:covering}
Given $P(n,k)$, find a polygonal trail whose union contains every grid point and whose number of segments is $h(n,k)$. The vertices may be arbitrary real points, including points outside $[0,n-1]^k$, and segment directions are unrestricted. The case $(n,k)=(3,2)$ is the nine-dot problem.
\end{problem}

Several subfamilies have already been solved. For $n\ge3$, the planar optimum $h(n,2)=2(n-1)$ is established in the 1970 treatment of Golomb and Selfridge \cite[Sections~2 and~4]{golomb}; it is also part of the 2014 classification of covering paths and trees for planar grids \cite[Theorem~3]{keszegh}. The constructions and lower bounds in \cite[Theorems~2.2--2.3]{rinaldi} give
\[
 h(2,k)=3\cdot2^{k-2}\qquad(k\ge2),
\]
with the same optimum attainable by closed trails. For the ternary grids, the clockwise-algorithm \cite{clockwise} attains $h(3,k)=(3^k-1)/2$. Subsection~\ref{sec:ternary} gives a short derivation of the matching lower bound from reduction modulo $3$. Thus, the hypercubic family includes both the vertices of a hypercube and the natural higher-dimensional version of the nine-dot puzzle.

Related covering problems impose additional geometric constraints. Rectilinear Hamiltonian tours in grids are studied in \cite{kranakis}, and uncrossing paths constrained to an axis-aligned bounding box in \cite{uncrossing}. For general planar point sets, covering paths lead to further extremal bounds and approximation questions \cite{planar,approx,biniaz,akitaya}. These planar extremal and approximation results do not determine the optimum for the hypercubic grids considered here, where long collinear groups and repeated visits to grid points are central to the minimum-link problem. Dumitrescu and T\'oth \cite{dumitrescutoth} study covering trees and distinguish the number of edges from the number of maximal collinear segments. They also consider connected arrangements of segments, which need not be traversable as a single covering trail with the same number of segments. Their restrictions on the directions of lines containing many grid points provide a precedent for Lemma~\ref{geo:lem:unit}. Their counting method also motivates the large-grid obstruction developed in Appendix~\ref{app:largegrid}, where the argument is adapted directly to polygonal trails allowing isolated self-intersections.

Our first result gives a geometric proof, together with a Lean formalization, of the hypercubic capacity bound stated in \cite[Theorem~2.1]{ripa}. The numerical bound and the principle of counting repetitions are prior results; the proof below supplies the geometric argument for every $n\ge3$ and every finite dimension $k\ge2$. The difficulty is that an external turn need not occur at a grid point, so consecutive segments need not repeat a grid point at their common endpoint. We address this by considering a counterexample with the least possible number of segments. Its segments would each contain at least $n-1$ grid points, and their grid-point sets would be pairwise disjoint. For $n\ge4$, an invariant of the external turns rules out this configuration; for $n=3$, a separate argument modulo $3$ gives the same capacity bound.

\begin{theorem}[Hypercubic capacity]\label{thm:maincapacity}
For integers $n\ge3$, $k\ge2$ and $r\ge1$, every polygonal trail of $r$ segments in $\mathbb R^k$ meets at most $1+r(n-1)$ distinct points of $P(n,k)$. Consequently,
\begin{equation}\label{eq:capacitylower}
 h(n,k)\ge\frac{n^k-1}{n-1}=1+n+\cdots+n^{k-1}.
\end{equation}
\end{theorem}
The proof occupies Section~\ref{sec:capacity}. The expression $n+(r-1)(n-1)$, obtained by counting the first segment separately as in \cite[Theorem~2.1]{ripa}, is the same capacity bound.

Our second result determines the optimum for $P(4,3)$, the case distinguishing the alternative general formulae proposed in \cite{conjecture}.
\begin{theorem}[The $4\times4\times4$ grid]\label{thm:main23}
The minimum-link covering number of $P(4,3)$ is $h(4,3)=23$. Furthermore, a trail of $22$ segments can cover exactly $63$ of its points, the maximum possible.
\end{theorem}

The earlier construction in \cite{small}, also discussed in \cite{conjecture}, already gives the upper bound $h_u(4,3)=23$. In October 2022, Alekseyev \cite{alekseyev} enumerated the $1492$ supporting lines and $13307$ intersection points determined by this grid and supplied Sage code. To obtain a complete lower-bound argument from a finite model, we must also justify why every hypothetical cover occurs in that model. Section~\ref{sec:four} proves that each segment of a $22$-segment cover would have two grid points visited by no other segment, thereby determining its support. We use exact rational incidence and order tests. Our reduction leaves $632$ initial pairs of states. The certified search excludes every continuation of those pairs, so no $22$-segment cover exists; together with the earlier construction, this proves the asserted optimum. The same section gives an explicit $22$-segment trail through $63$ points, while the companion technical note \cite{technical} supplies the full computational account.

\begin{remark}[Allowing retracing]\label{rem:walks}
The capacity bound in Theorem~\ref{thm:maincapacity} and the lower bounds in Theorem~\ref{thm:main23} and Proposition~\ref{prop:five32} remain valid for walks allowing retracing. Their proofs apply to the more general polygonal chains defined in Section~\ref{sec:capacity}, which also allow degenerate segments, while the constructions attaining the upper bounds used here are trails.
\end{remark}

Our formalization uses Lean 4 \cite{lean} and mathlib \cite{mathlib}. The geometric reduction and the verification of the remaining cases have distinct roles: here, the capacity proof is checked by the kernel, whereas the numerical exclusions use native reflection \cite{native}. Section~\ref{sec:formal} establishes the lower bound $h(5,3)\ge32$, presents an optimal $23$-segment trail in $[0,4]^3$, and describes the trust assumptions of the formal verifications. It then revisits the proposed general formulae in light of the large-grid obstruction proved in Appendix~\ref{app:largegrid}.

Table~\ref{tab:bounds} places the two numerical improvements among the known bounds for small hypercubic grids; the exact formulae above continue its first two rows to every $k\ge2$.

\begin{table}[H]
\centering
\renewcommand{\arraystretch}{1.35}
\setlength{\tabcolsep}{13pt}
\begin{tabular}{|c|cccc|}
\hline
\diagbox[width=6em,height=2.8em]{$n$}{$k$} & $2$ & $3$ & $4$ & $5$\\
\hline
$2$ & \cellcolor{exactcell}$3$ & \cellcolor{exactcell}$6$ & \cellcolor{exactcell}$12$ & \cellcolor{exactcell}$24$\\
$3$ & \cellcolor{exactcell}$4$ & \cellcolor{exactcell}$13$ & \cellcolor{exactcell}$40$ & \cellcolor{exactcell}$121$\\
$4$ & \cellcolor{exactcell}$6$ & \cellcolor{exactcell}$\mathbf{23}$ & \cellcolor{boundcell}$[85,95]$ & \cellcolor{boundcell}$[341,383]$\\
$5$ & \cellcolor{exactcell}$8$ & \cellcolor{boundcell}$[\mathbf{32},36]$ & \cellcolor{boundcell}$[156,184]$ & \cellcolor{boundcell}$[781,924]$\\
\hline
\end{tabular}
\caption{Bounds for $h(n,k)$, adapted from \cite[Table~1]{conjecture}. Green cells give exact values and orange cells give intervals $[h_l(n,k),h_u(n,k)]$. Bold numerals mark the two improvements established here: the exact value $h(4,3)=23$, replacing $[21,23]$, and the lower bound $h(5,3)\ge32$, replacing $31$. The interval $[32,36]$ does not assert an exact value.}\label{tab:bounds}
\end{table}

For $n\ge3$ and $k\ge3$, the construction in \cite{conjecture} and Theorem~\ref{thm:maincapacity} give
\begin{equation}\label{eq:boundsallk}
 \frac{n^k-1}{n-1}\le h(n,k)\le\bigl(h_u(n,3)+1\bigr)n^{k-3}-1.
\end{equation}
In particular, $h_u(4,3)=23$ and $h_u(5,3)=36$ yield the displayed upper bounds in dimensions $4$ and $5$, while the corresponding lower bounds follow from~\eqref{eq:capacitylower}.

\FloatBarrier

\section{A geometric proof of the hypercubic capacity bound}\label{sec:capacity} This section establishes the capacity bound in Theorem~\ref{thm:maincapacity}, together with its extension to chains allowing retracing and degenerate segments, as recorded in Remark~\ref{rem:walks}. For this purpose, a \emph{polygonal chain} of $r\ge1$ segments is given by vertices $v_0,\ldots,v_r\in\mathbb R^k$ and the consecutive closed segments $L_i=[v_{i-1},v_i]$, with retracing and degenerate segments allowed. Every trail is therefore a chain. Working in this larger class lets us apply the same capacity estimate to any contiguous part left by cutting a chain, and also allows us to extend shorter covers by repeating their final vertex. 

\begin{definition}[Grid-point terminology]\label{def:gridterminology}
Let $n\ge3$, $k\ge2$, and let $L_1,\ldots,L_r$ be a polygonal chain. Write
\[
 S_i=L_i\cap P(n,k),\qquad q_i=\card{S_i}\qquad(1\le i\le r).
\]
A segment is \emph{strong} if it contains at least two grid points and \emph{weak} otherwise. It is \emph{rich} if it contains at least $n-1$ grid points and \emph{full} if it contains $n$.

Relative to a specified family of segments, a grid point is \emph{private} to one segment if it belongs to that segment and to no other member of the family. A chain of $r$ segments is \emph{saturated} if it covers $1+r(n-1)$ grid points.
\end{definition}

The term \emph{full} agrees with Dumitrescu and T\'oth's terminology for lines containing $n$ grid points \cite[Section~3]{dumitrescutoth}. Their \emph{heavy} lines contain more than $\lceil n/2\rceil$ grid points, whereas our \emph{rich} segments have the specific threshold $n-1$.

Throughout the proofs, $i,j$ index segment positions, while $\alpha,\beta,\gamma$ index coordinates in $\{1,\ldots,k\}$. Since a nonconstant coordinate is injective on a line, every segment contains at most $n$ grid points; a degenerate segment contains at most one. Summing these individual bounds gives $rn$, so the desired estimate requires a further argument about the way the segments are joined. The following counting lemma isolates the only configuration that could violate the bound: its least-counterexample argument is the strong induction used in the formalization.

\begin{lemma}[Properties of a least counterexample]\label{geo:lem:minimal}
Fix $n\ge3$ and $k\ge2$. If a chain with the least possible number $r$ of segments violates the bound $1+r(n-1)$, then every segment contains at least $n-1$ grid points and the sets $S_1,\ldots,S_r$ are pairwise disjoint.
\end{lemma}
\begin{proof}
By minimality, every proper contiguous subchain of $\ell$ segments covers at most $1+\ell(n-1)$ points; this inequality is also valid for an empty subchain. If $q_i\le n-2$, we may apply it to the portions before and after $L_i$. Together with $L_i$, these portions cover at most
\[
1+(i-1)(n-1)+(n-2)+1+(r-i)(n-1)=1+r(n-1)
\]
points, a contradiction. Thus, every segment is rich.

Similarly, if $S_i\cap S_j\ne\varnothing$ for some $i<j$, splitting immediately before $L_j$ gives two proper subchains with a common grid point. Applying the bound to each and subtracting this repetition gives
\[
[1+(j-1)(n-1)]+[1+(r-j+1)(n-1)]-1=1+r(n-1),
\]
again contradicting the assumed coverage. Therefore, the grid-point sets are pairwise disjoint.
\end{proof}

\subsection{External turns for \texorpdfstring{$n\ge4$}{n at least 4}}
Throughout this subsection, assume $n\ge4$. The richness condition then severely restricts the geometry of a least counterexample. We first determine the possible directions, then locate the turns between consecutive segments. The direction restriction is a special case of \mbox{\cite[Observation~2(3)]{dumitrescutoth}} since $n-1>\lceil n/2\rceil$ for $n\ge4$; we include a direct proof for completeness. These local restrictions will force a repeated grid point between two full segments, even when those segments are not consecutive.

\begin{lemma}[Unit directions]\label{geo:lem:unit}
The grid points of a rich segment form a consecutive interval in an integer direction belonging to $\{-1,0,1\}^k\setminus\{0\}$.
\end{lemma}
\begin{proof}
Order the grid points by a nonconstant coordinate along the line. Among at least $n-1$ distinct integer values in $\{0,\ldots,n-1\}$, two consecutive values differ by one, so the difference $d$ of the corresponding grid points gives a direction in which every other grid point on the line has integer parameter. There are at least $n-1$ such parameters, and consequently each nonzero coordinate of $d$ satisfies $(n-2)|d_{\alpha}|\le n-1$. Since $n\ge4$, this forces $|d_{\alpha}|=1$. Convexity of both the cube and the segment then shows that the parameters form a consecutive interval.
\end{proof}

The next lemma starts with an arbitrary real turning point. Its integrality is a consequence of richness and disjointness, rather than a restriction imposed on the search space.
\begin{lemma}[Disjoint rich segments]\label{geo:lem:turn}
Two consecutive rich segments with disjoint grid-point sets have distinct supporting lines. Their common endpoint is an integer point outside $[0,n-1]^k$ and inside $[-1,n]^k$. On either segment, the first grid point from that endpoint is one unit-direction step away.
\end{lemma}
\begin{proof}
A supporting line contains at most $n$ grid points, so it cannot contain two disjoint sets of at least $n-1$ grid points. The two supports are therefore distinct. Orient their unit directions $d,e$ away from the common endpoint $v$. In any coordinate active in both directions, the signs must agree: otherwise, the two sets of at least $n-1$ distinct coordinate values would lie on opposite sides of $v_{\alpha}$ and have at most one value in common. This is impossible since two subsets of $\{0,\ldots,n-1\}$, containing at least $n-1$ elements each, must share at least $n-2\ge2$ elements.

The sets of active coordinates cannot be equal, for the sign agreement would then give $d=e$ and coincident supports. 
Thus, choose a coordinate $\alpha$ active on one line and constant on the other; its value at $v$ is an integer in $\{0,\ldots,n-1\}$ since the other line contains grid points. If $p=v+tu$ is a grid point on the active line, where $u$ is its oriented unit direction, then
\[
t=\frac{p_\alpha-v_\alpha}{u_\alpha}\in\mathbb Z,
\]
since $u_\alpha=\pm1$. Consequently, $v=p-tu\in\mathbb Z^k$. Moreover, $v$ lies outside the cube, since otherwise it would be a grid point common to the two segments.

On the line with the exclusive active coordinate, the grid points have positive integer parameters measured from $v$. The first parameter must be one: since $v_\alpha\in\{0,\ldots,n-1\}$ and $|u_\alpha|=1$, starting at parameter two would allow at most $n-2$ grid points. Thus, one unit-direction step reaches the first grid point and, in particular, $v\in[-1,n]^k$. To obtain the same conclusion on the other segment, choose a coordinate in which $v$ lies outside the cube. Its value is $-1$ or $n$, so a unit step toward the grid enters the allowed range. The step also lies in range in every remaining coordinate: those already inside remain between their values at $v$ and at a later grid point, while any other external coordinate is handled in the same way. Hence, the first step on this segment is a grid point as well.
\end{proof}

The displacement between successive external turns now gives the invariant needed for the capacity bound.
\begin{lemma}[The external-turn invariant]\label{geo:lem:anchor}
Consider a contiguous part of a chain beginning and ending with two distinct full segments, with only $(n-1)$-point segments in between. If consecutive segments have disjoint grid-point sets, the two full segments share a grid point.
\end{lemma}
\begin{proof}
By Lemma~\ref{geo:lem:turn}, each intervening turn is an external integer point one step beyond the nearest grid point on either incident segment. An intermediate segment therefore has $n-1$ grid points and one additional step at each end, so its endpoints differ by $nd$ for its oriented unit direction $d$. It follows that all the intervening turning points have the same coordinatewise residue modulo $n$. Let $p_*\in P(n,k)$ be the representative of that residue.

At the external endpoint $v$ of either full segment, orient the unit direction $u$ toward the grid. The $n$ grid points on that segment are $v+u,\ldots,v+nu$; the last is congruent to $v$ modulo $n$ and hence equals $p_*$. Thus, both full segments contain $p_*$.
\end{proof}

\begin{remark}[A nonadjacent repetition]
For $n=4$, consider the planar chain
\[
(3,0)\longrightarrow(-1,0)\longrightarrow(3,4)\longrightarrow(3,0).
\]
Its successive segments contain $4,3,4$ grid points. Although neither turn repeats a grid point, the first and last segments share $(3,0)$, the common residue of the turning points modulo $4$. Their union therefore contains ten points, exactly $1+3(4-1)$. This illustrates why the invariant concerns repetitions along a whole part of the chain, not just at its consecutive turns.
\end{remark}

\begin{theorem}\label{geo:thm:cubic}
For $n\ge4$ and $k\ge2$, every chain of $r\ge1$ segments in $\mathbb R^k$ meets at most $1+r(n-1)$ points of $P(n,k)$.
\end{theorem}
\begin{proof}
Suppose a counterexample exists and choose one with the least number of segments. By Lemma~\ref{geo:lem:minimal}, its segments are rich and their grid-point sets are pairwise disjoint. If at most one segment is full, their cardinalities sum to at most $r(n-1)+1$, which is already a contradiction. Otherwise, take two successive full segments in the ordered list of full segments. Every segment between them contains exactly $n-1$ grid points, so Lemma~\ref{geo:lem:anchor} gives a grid point common to the two full segments. This contradicts disjointness and completes the proof.
\end{proof}

\subsection{The ternary grids: a finite-field argument}\label{sec:ternary}
When $n=3$, two grid points do not force a unit direction; accordingly, we use a different description of the segments. Reduction modulo $3$ identifies $P(3,k)$ bijectively with $\mathbb F_3^k$. Denoting coordinatewise reduction by a bar, two distinct grid points $a,a+d$ determine a three-point affine line $\bar a+\mathbb F_3\bar d$. Here $d_{\alpha}\in\{-2,-1,0,1,2\}$ and $\bar d\ne0$, with no primitive-direction assumption. For any integer point $x$ on the real support, eliminating the real line parameter gives
\[
d_\alpha(x_\beta-a_\beta)
=d_\beta(x_\alpha-a_\alpha)
\qquad(1\le\beta\le k).
\]
Reducing modulo $3$ and choosing $\alpha$ with $\bar d_\alpha\ne0$ shows that $\bar x$ belongs to $\bar a+\mathbb F_3\bar d$. Since two distinct points determine a line over $\mathbb F_3$, the completed line is independent of the chosen pair of grid points.

The grid representative of a point added by completion need not lie on the original segment or even on its real supporting line. Passing to the completed line nevertheless only enlarges the reduced set to be bounded. Its useful property is that consecutive completed lines always intersect.

\begin{lemma}[Completed ternary lines intersect]\label{geo:lem:ternary}
If two consecutive real segments each contain at least two points of $P(3,k)$, their completed affine lines over $\mathbb F_3$ intersect.
\end{lemma}
\begin{proof}
Let $v$ be the common endpoint, and on each segment choose two grid points in order away from $v$, denoted by $a,a+d$ and $b,b+e$. Thus
\[
 v=a-\lambda d=b-\mu e,\qquad \lambda,\mu\ge0.
\]
If some minor $D=d_{\alpha}e_{\beta}-d_{\beta}e_{\alpha}$ is not divisible by $3$, Cramer's rule applied to $a-b=\lambda d-\mu e$ shows that $D\lambda$ and $D\mu$ are integers. For every coordinate $\gamma$, the identity
\[
D(a_{\gamma}-b_{\gamma})
=(D\lambda)d_{\gamma}-(D\mu)e_{\gamma}
\]
can therefore be reduced modulo $3$. Since $\bar D\ne0$, division by $\bar D$ gives parameters in $\mathbb F_3$ defining a common point of the two completed lines.

Otherwise, every minor is divisible by $3$. If at least one is nonzero, its value must be one of $\pm3,\pm6$, since the direction entries belong to $\{-2,-1,0,1,2\}$. These possibilities give $80$ cases: $24$ each with determinant $3$ and $-3$, and $16$ each with determinant $6$ and $-6$. In every case, the equations
\[
 a_{\alpha}-\lambda d_{\alpha}=b_{\alpha}-\mu e_{\alpha},\qquad
 a_{\beta}-\lambda d_{\beta}=b_{\beta}-\mu e_{\beta},
\]
together with $\lambda,\mu\ge0$ and the coordinate bounds
\[
 0\le a_{\gamma},\ a_{\gamma}+d_{\gamma},\ b_{\gamma},\ b_{\gamma}+e_{\gamma}\le2
 \qquad(\gamma=\alpha,\beta),
\]
force $\lambda=\mu=0$. Up to exchanging the two coordinates, reflecting either coordinate within $\{0,1,2\}$, and interchanging the two segments, these $80$ direction configurations reduce to six classes. Four of these classes admit no feasible configuration with $\lambda,\mu\ge0$, while the remaining two force $\lambda=\mu=0$. Appendix~\ref{app:ternarycases} gives all six representatives and elementary coordinate inequalities for each case. This implication is checked in the Lean lemma \path{ternary_bad_determinant} by integer case splitting and linear arithmetic; throughout the check, the parameters $\lambda,\mu$ remain real. Hence $a=b=v$, which gives a common point after reduction.

It remains to consider the case in which all minors vanish. The vectors $d,e$ are then proportional, and their real supporting lines coincide since both contain $v$. Choose $\alpha$ with $\bar d_{\alpha}\ne0$. Collinearity gives
\[
 (b_{\beta}-a_{\beta})d_{\alpha}=(b_{\alpha}-a_{\alpha})d_{\beta}
 \qquad\text{for every }\beta.
\]
After reduction, division by $\bar d_{\alpha}$ shows that $\bar b$ belongs to the completed line through $a$ as well as to its own completed line. This also covers collinear segments and arbitrary real turning points.
\end{proof}

\begin{theorem}\label{geo:thm:ternary}
For every $k\ge2$, a chain of $r\ge1$ segments meets at most $2r+1$ points of $P(3,k)$.
\end{theorem}
\begin{proof}
First suppose that every segment contains at least two grid points and complete each segment to its three-point affine line over $\mathbb F_3$. By Lemma~\ref{geo:lem:ternary}, each line after the first meets its predecessor and therefore adds at most two points to the union. The union of the completed lines has at most $3+2(r-1)=2r+1$ points and contains the bijective image of every visited grid point, proving the bound in this case.

For a general chain, suppose that the bound fails and take a counterexample with the least number of segments. Lemma~\ref{geo:lem:minimal}, applied with $n=3$, forces every segment to contain at least two grid points, so the preceding argument excludes this counterexample as well.
\end{proof}

\begin{proof}[Proof of Theorem~\ref{thm:maincapacity}]
Theorems~\ref{geo:thm:cubic} and~\ref{geo:thm:ternary} establish the capacity bound for chains and hence for trails. A covering chain must therefore satisfy $n^k\le1+r(n-1)$, which gives~\eqref{eq:capacitylower}.
\end{proof}

In the ternary case, the clockwise-algorithm attains this lower bound and gives the following exact value in every dimension.
\begin{corollary}[The ternary optimum]\label{cor:ternary}
For every $k\ge2$, $h(3,k)=(3^k-1)/2$.
\end{corollary}
\begin{proof}
Theorem~\ref{geo:thm:ternary} requires $3^k\le2r+1$ for any cover, while the clockwise-algorithm~\cite{clockwise} constructs covering trails with exactly $(3^k-1)/2$ segments for every $k\ge2$.
\end{proof}

\section{The optimum for \texorpdfstring{$P(4,3)$}{P(4,3)}}\label{sec:four}
For $P(4,3)$, the general capacity bound gives only $h(4,3)\ge21$. To reach $23$, we refine the capacities of shorter chains and use them to show that every segment of a hypothetical $22$-segment cover has two private grid points. This property makes the final enumeration complete since it determines every supporting line from points of the grid.

Accordingly, throughout this section, let $C(r)$ denote the largest number of distinct points of $P(4,3)$ covered by a chain of $r$ segments, with $C(0)=0$ for the empty segment family. In the terminology of Definition~\ref{def:gridterminology}, a saturated chain of $r$ segments covers $3r+1$ points.

\subsection{Capacity refinements and disconnected packings}
Deleting an internal segment of a polygonal chain leaves two chains, and further deletions may leave more components.
In particular, deleting all weak segments leaves components whose segments are all strong in the sense of Definition~\ref{def:gridterminology}; we call any such chain a strong chain. We therefore need finite exclusions for families of strong chains as well as bounds for a single chain. The components in such a family need not be joined to one another; disjointness, when required, refers to their covered grid-point sets.

\begin{lemma}[Finite exclusions]\label{lem:finite}
The exact finite checks for $P(4,3)$ establish the following statements:
\begin{enumerate}
\item A strong chain of $20$ segments covers at most $59$ points and, assuming $C(20)\le59$, a strong chain of $21$ segments covers at most $62$ points;
\item There are no two disjoint saturated strong chains of total length $19$ or $20$, and no four disjoint saturated strong chains of total length $19$;
\item Three nonempty strong chains of total length $20$ cannot cover at least $62$ points;
\item Two nonempty strong chains of total length $21$ cannot cover at least $63$ points if every segment has a point private to it in the two-chain family.
\end{enumerate}
\end{lemma}
\begin{proof}[Proof (finite verification)]
Encode a strong segment by its first and last grid points in the direction of traversal. Since these points are distinct, their ordered pair determines the supporting line and the exact set of grid points visited by the segment, namely those lying between them, inclusive. This gives $4032=64\cdot63$ oriented states, corresponding to the ordered pairs of distinct grid points. For a transition between two states with distinct supports, the corresponding lines must meet at or beyond the last grid point of the first state and at or before the first grid point of the second. Extending the segments to that turn must also preserve their encoded grid-point sets exactly. Since each supporting line is determined by two integer grid points, the intersection of any two distinct such lines, when it exists, has rational coordinates. Thus, all incidence, order and set-equality conditions can be tested using exact rational arithmetic; the preliminary adjacency table is proved to contain every geometrically admissible transition.

To justify this encoding, we must exclude consecutive collinear segments in the families under consideration. For a family of $m$ nonempty strong chains with $r$ segments in total, the sum of their capacity bounds is $3r+m$. If their union covers $U$ distinct grid points, the deficit $3r+m-U$ consists of the shortfalls from the individual capacity bounds together with repetitions between components.
This deficit is zero for the configurations in item~(2) and at most one for those in item~(3). Consequently, every component of $\ell$ segments covers at least $3\ell$ points.

Consecutive collinear segments can be treated in two cases. If they form a forward continuation, they can be merged without changing the covered set, leaving a chain of $\ell-1$ segments with capacity at most $3(\ell-1)+1=3\ell-2$. Otherwise, they retrace part of the same supporting line, and one segment is contained in the other. Deleting the contained segment preserves coverage and leaves at most two nonempty chains with $\ell-1$ segments in total, whose combined capacity is at most $3(\ell-1)+2=3\ell-1$. Both bounds contradict the required coverage of at least $3\ell$ points.

The same argument applies to the hypothetical counterexamples in item~(1): a $20$-segment chain covering at least $60$ points or a $21$-segment chain covering at least $63$ points would fall short of its general capacity bound by at most one. In item~(4), the difference is at most two, so each component covers at least $3\ell-1$ points: this still excludes merging, and the explicit privacy hypothesis excludes retracing. Consequently, every relevant real family admits the finite encoding.

The $20$-segment exclusion has $116$ initial states up to coordinate permutations and reflections. Once $C(20)\le59$ has been established for arbitrary real chains, deleting either end segment of a $21$-segment chain covering at least $63$ points leaves a $20$-segment chain covering at most $59$ points. Each end segment must therefore contain at least $63-59=4$ private grid points. This reduces the initial states for the second exclusion to six symmetry representatives.

For the five packing exclusions in items~(2)--(4), the finite search also considers the possible distributions of the total number of segments among the nonempty components, subject to the capacity restrictions above. Admissible transitions are checked within each component, together with the disjointness or privacy conditions required by the corresponding item.

Lean verifies the completeness of these geometric reductions and the pruning rules, while native reflection checks the remaining finite predicates. Their exact definitions and certificates are given in \cite{technical}.
\end{proof}

The first application extends the single-chain bounds to chains containing weak segments. The order matters: the $20$-segment bound must be proved before it can justify the six initial states for $21$ segments.
\begin{lemma}\label{lem:weakcapacity}
The refined capacities satisfy $C(20)\le59$ and $C(21)\le62$.
\end{lemma}
\begin{proof}
Let $r\in\{20,21\}$ and suppose that a chain of $r$ segments covers at least $3r$ points. If $s$ segments are weak, deleting them leaves at most $s+1$ strong components. Applying Theorem~\ref{geo:thm:cubic} to these components and allowing one point for each deleted segment bounds the original coverage by $3(r-s)+(s+1)+s=3r+1-s$.
It follows that $s\le1$.

If the unique weak segment occurs at an end, the remaining $r-1$ segments form one chain, so the total coverage is at most $3(r-1)+1+1=3r-1$, a contradiction.

Hence, the weak segment must be internal; deleting it leaves two nonempty strong chains of total length $r-1$, whose combined capacity is at most $3(r-1)+2=3r-1$.

To reach at least $3r$ points after restoring the weak segment, both chains must be saturated, their grid-point sets must be disjoint, and the weak segment must contribute one additional grid point. Lemma~\ref{lem:finite}(2) excludes that configuration, so all segments are strong.

Item~(1) now excludes $r=20$, proving $C(20)\le59$; with this premise available, it excludes $r=21$ as well.
\end{proof}

We can therefore apply the packing exclusions to a hypothetical cover. The following lemma is the essential link between the continuous geometric problem and the final finite search.
\begin{lemma}\label{lem:twoprivate}
Every segment of a $22$-segment cover of $P(4,3)$ would have at least two private grid points. Consecutive supporting lines would be distinct.
\end{lemma}
\begin{proof}
We first eliminate weak segments. If their number is $s$, deletion and the capacity bound give $64\le3(22-s)+(s+1)+s=67-s$, so $s\le3$. For $s=3$, equality forces four disjoint saturated strong components of total length $19$, contrary to Lemma~\ref{lem:finite}(2). For $s=2$, the remaining $20$ strong segments must cover at least $62$ points, since the two weak segments contribute at most two grid points. If these strong segments form a single component, its capacity is at most $3\cdot20+1=61$. If they form two nonempty components, their combined capacity is $3\cdot20+2=62$, so equality requires both components to be saturated and their grid-point sets to be disjoint, contrary to Lemma~\ref{lem:finite}(2). Three nonempty components are excluded by Lemma~\ref{lem:finite}(3). Thus, at most one segment is weak.

Suppose that there is exactly one weak segment. Every other segment must have a private point: otherwise, deleting a strong segment without private points as well as the weak segment would leave $20$ strong segments in at most three components covering at least $63$ points. One or two components cannot cover that many points, and Lemma~\ref{lem:finite}(3) excludes three. The weak segment cannot occur at an end since the other $21$ segments cover at most $62$ points by Lemma~\ref{lem:weakcapacity}. Deleting the internal weak segment therefore leaves two nonempty strong components of total length $21$ covering at least $63$ points. Each remaining segment retains a private point in this two-component family, contrary to Lemma~\ref{lem:finite}(4). Hence, every segment of the cover is strong.

Next suppose that some segment has no private point. It must again be internal, and its deletion leaves two strong components covering all $64$ points. Deleting any further segment from this family leaves at most three components and $20$ segments, whose total capacity is at most $63$. Every segment of the two-component family must consequently have a private point, so item~(4) gives another contradiction.
Thus, we have proved that every original segment has at least one private point.

It remains to exclude the possibility that a segment has exactly one private point. An end segment must contain at least two private grid points since deleting it leaves a $21$-segment chain covering at most $62$ of the original $64$ grid points. If an internal segment had exactly one private point, its deletion would leave two nonempty strong components covering at least $63$ points, and every remaining segment would retain a private point. Item~(4) excludes this final case, proving the required two-point privacy.

Finally, consecutive collinear segments either retrace, in which case one is contained in the other and has no private point, or merge into a forward continuation. The latter would give a $21$-segment cover, contradicting $C(21)\le62$.

Therefore, consecutive supports are distinct.
\end{proof}

\subsection{Completeness of the final enumeration}\label{sec:enumeration}
Let $M_\sigma$ denote the grid-point set encoded by a state $\sigma$. By Lemma~\ref{lem:twoprivate}, every transition $\sigma\to\tau$ of a hypothetical cover satisfies
\[
 |M_\sigma\setminus M_\tau|\ge2,\qquad |M_\tau\setminus M_\sigma|\ge2.
\]
Adding these necessary conditions to the transition table gives a graph with $54624$ arcs. Moreover, deleting the first two segments leaves a $20$-segment chain covering at most $59$ points, so the first pair has at least five points private to it and, in particular, at least five points in its union. Normalization under the $48$ cube symmetries reduces these initial pairs to the certified list of $632$ two-state prefixes. No box or lattice restriction is imposed on the turning points: the internal turns are rational because they are intersections of the determined supports.

The remaining pruning conditions express how much coverage a prefix must still make possible. The guarded search uses the generic capacity upper bound $c$, defined by $c(0)=0$ and $c(r)=3r+1$ for $r\ge1$. In a prefix of $\ell$ segments, every segment must retain at least two points private within that prefix.

More generally, consider a block consisting of segments $i+1,\ldots,j$, where $0\le i\le j\le\ell$. In any completed $22$-segment cover, the first $i$ segments can cover at most $c(i)$ grid points, while the last $22-j$ segments can cover at most $c(22-j)$ points. Since the entire grid contains $64$ points, at least $64-c(i)-c(22-j)$ points must be covered by the selected block and by no segment outside it. In particular, the block must already contribute at least this many points outside the other segments recorded in the current prefix.

This requirement includes the prefix bound $64-c(22-\ell)$ when $i=0$ and $j=\ell$. Since a block's private contribution can only decrease as the prefix is extended, failure of any of these conditions excludes a completed cover. The sharper bounds for $20$ and $21$ segments enter the preceding geometric and initial-state reductions; this guard uses the generic function $c$.

Here, the budget is the allowed capacity deficit $3\cdot22+1-64=3$, not a limit on runtime or on the number of visited search nodes. The remaining depth counts segments still to append; a complete leaf is accepted only if all $64$ points are covered and its guards hold. An interrupted computation supplies no negative certificate.

For each of the $632$ roots, native evaluation of the budgeted search finds no admissible completion. The formal verification proves that any accepting continuation of the guarded search would also be accepted by the budgeted search and that the machine-word implementation agrees with the formal bit-vector computation. These implications transfer the computed verdicts to the guarded root predicates and hence to the original geometric problem.

\begin{figure}[!htbp]
\centering
\includegraphics[width=\linewidth]{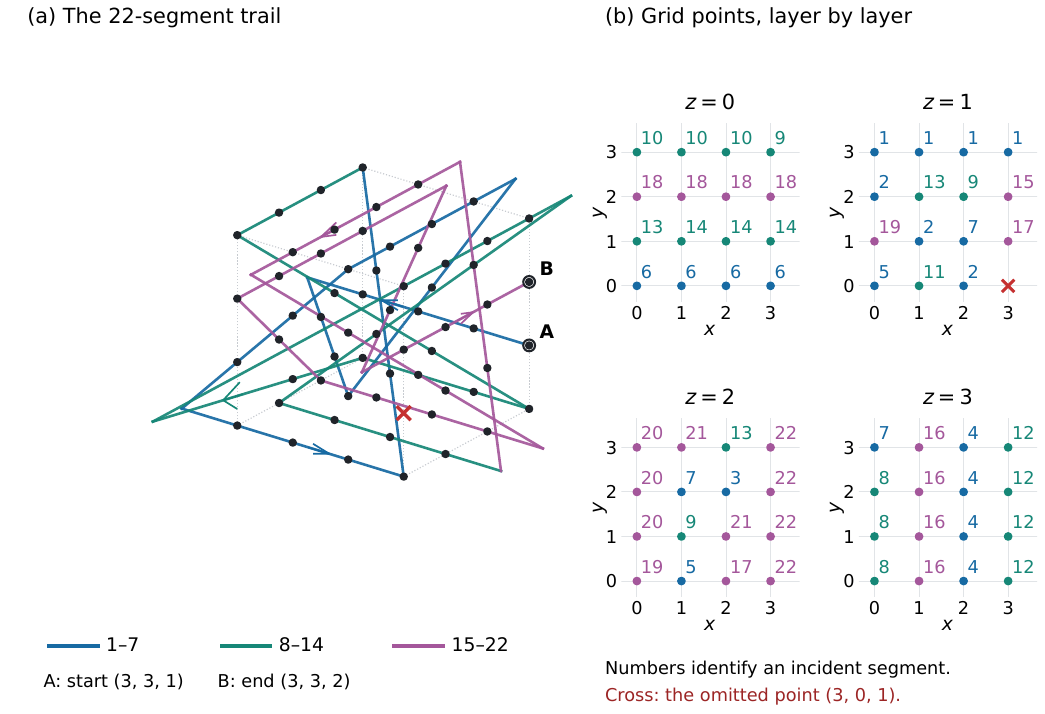}
\captionsetup{singlelinecheck=false}
\caption[An explicit $22$-segment trail through $63$ of the $64$ grid points]{An explicit $22$-segment trail through $63$ points of $P(4,3)$. Its ordered vertices are:
\[
\begin{aligned}
&(3,3,1)\to(-1,3,1)\to(2,0,1)\to(2,4,3)\to(2,0,3)\to(-1,0,0)\\
 &\to(3,0,0)\to(0,3,3)\to(0,0,3)\to(3,3,0)\to(0,3,0)\to(3,-6,3)\\
 &\to(3,4,3)\to(0,1,0)\to(4,1,0)\to(1,4,3)\to(1,-1,3)\to(4,2,0)\\
 &\to(0,2,0)\to(0,0,2)\to(0,5,2)\to(3,-1,2)\to(3,3,2).
\end{aligned}
\]
Panel (a) shows the full trail, with start $A=(3,3,1)$ and end $B=(3,3,2)$. Panel (b) displays the four horizontal grid layers: each small number identifies a segment passing through the corresponding grid point, with segments numbered from $1$ in the displayed sequence. The only omitted point is $(3,0,1)$, marked by a cross. Appending that point to the sequence gives a $23$-segment covering trail. These are the coordinates used in the Lean witness; all displayed incidences were also recomputed using exact rational arithmetic.} \label{fig:trail22}
\end{figure}

\begin{proof}[Proof of Theorem~\ref{thm:main23}]
A $22$-segment cover would satisfy Lemma~\ref{lem:twoprivate} and therefore occur in the complete enumeration, contradicting the finite certificate. A shorter covering chain could be extended to $22$ segments by repeating its final vertex, so it is excluded as well. The earlier construction \cite{small,conjecture} attains $23$ segments and establishes the matching upper bound.

For the partial-cover statement, Figure~\ref{fig:trail22} displays the explicit witness retained in the Lean artifact: its first $22$ segments cover exactly $63$ points, and appending the omitted point gives a $23$-segment covering trail. Since no $22$-segment chain covers all $64$ points, this partial coverage is maximal.
\end{proof}

Thus, the refined capacities established in this section are $C(20)\le59$, $C(21)\le62$, and $C(22)=63$.

\section{Further results and formal verification}\label{sec:formal}
The preceding results have two further consequences: for $P(5,3)$, the capacity theorem reduces the exclusion of $31$ segments to a separate finite computation; for $P(4,3)$, the unrestricted lower bound certifies the optimality of a covering trail contained in the enlarged box $[0,4]^3$.

\subsection{A lower bound for \texorpdfstring{$P(5,3)$}{P(5,3)}}\label{sec:five}
Unlike the argument for $P(4,3)$, the next reduction does not require the disconnected-packing computations of Section~\ref{sec:four}. The general capacity bound alone supplies the privacy and prefix conditions needed for the finite exclusion.
\begin{proposition}\label{prop:five32}
Every polygonal trail covering $P(5,3)$ has at least $32$ segments. In particular, $h(5,3)\ge32$.
\end{proposition}
\begin{proof}
Theorem~\ref{geo:thm:cubic} excludes $30$ or fewer segments, so suppose that a chain of $31$ segments covers all $125$ grid points. Removing one segment leaves at most two chains with total capacity at most $4\cdot30+2=122$; hence, every segment has at least three private grid points. Similarly, removing two consecutive segments leaves at most two chains of total length $29$, whose combined capacity is at most $4\cdot29+2=118$. The removed pair must therefore cover at least $125-118=7$ points not visited by the remaining segments. Moreover, for every integer $\ell$ with $1\le\ell\le31$, the suffix after the first $\ell$ segments has capacity at most $4(31-\ell)+1$. Consequently, the first $\ell$ segments must cover at least $4\ell$ grid points. Consecutive supports cannot coincide: retracing would contradict privacy, while merging a forward continuation would produce a $30$-segment cover.

These properties allow us to encode each segment by its first and last grid points, using the exact-set, incidence and order conditions from Section~\ref{sec:four}.

There are $2676$ oriented states containing at least three grid points. Requiring each of two consecutive states to contain at least three grid points not covered by the other, and their union to contain at least seven grid points, leaves $12264$ admissible transitions. By applying coordinate reflections and permutations, we may assume that the coordinates of the first grid point are nondecreasing and at most $2$. Since the prefix bound requires the first segment to contain at least four grid points, only $55$ initial states remain.

The search proceeds through prefixes, recording the last state and the union of visited points. At depth $\ell$ it retains only unions of size at least $4\ell$ and merges identical pairs of last state and union. No pair survives to depth $21$. Lean verifies the state masks, the completeness of the transition and initial-state tables, and the fact that merging identical pairs preserves every admissible continuation. Native reflection checks that the search returns no surviving states, thereby excluding the existence of a $31$-segment cover.
\end{proof}

\subsection{An optimal trail in the enlarged box}\label{sec:boxed}

Beyond the unrestricted covering problem, it is natural to consider variants that impose additional geometric constraints on the covering trail.

For the $2\times2\times\cdots\times2$ hypercubic grids, the minimum number of segments can be attained even by covering circuits (closed polygonal trails) that visit every grid point exactly once and never revisit any vertex of the polygonal trail except when returning to the starting point to close the circuit \cite[Theorem~2.3 and Corollary~2.4]{rinaldi}. Thus, these additional restrictions do not increase the optimum for this family. The stronger requirement of having no self-intersections can also be imposed without increasing this optimum: the case $k=2$ is elementary, the original MathOverflow question presents explicit simple closed coverings for $k=3,4$ \cite{mo-simple}, and a subsequent answer extends the construction to every $k\ge5$ \cite{mo-simple-answer}.

The following result provides another example, this time involving a restriction on the region containing the trail.

The OEIS sequence A225227 \cite{oeis} considers a constrained variant of the same covering problem: a trail through $P(n,3)$ must stay in $[0,n]^3$, a box of side length $n$ rather than $n-1$. The grid remains $\{0,\ldots,n-1\}^3$, and the vertices of the trail may still have arbitrary real coordinates. For $n=4$, the construction below attains the unrestricted optimum within this enlarged box, as the clockwise-algorithm does for $n=3$ \cite{clockwise}. It is also recorded, with its certificate, in~\mbox{\cite[Section~3.2]{technical}} and in the sequence entry.

\begin{proposition}[The box-constrained optimum]\label{prop:boxed}
The minimum number of segments of a trail covering $P(4,3)$ and contained in $[0,4]^3$ is $23$.
\end{proposition}

\begin{figure}[!htbp]
\centering
\includegraphics[width=\linewidth]{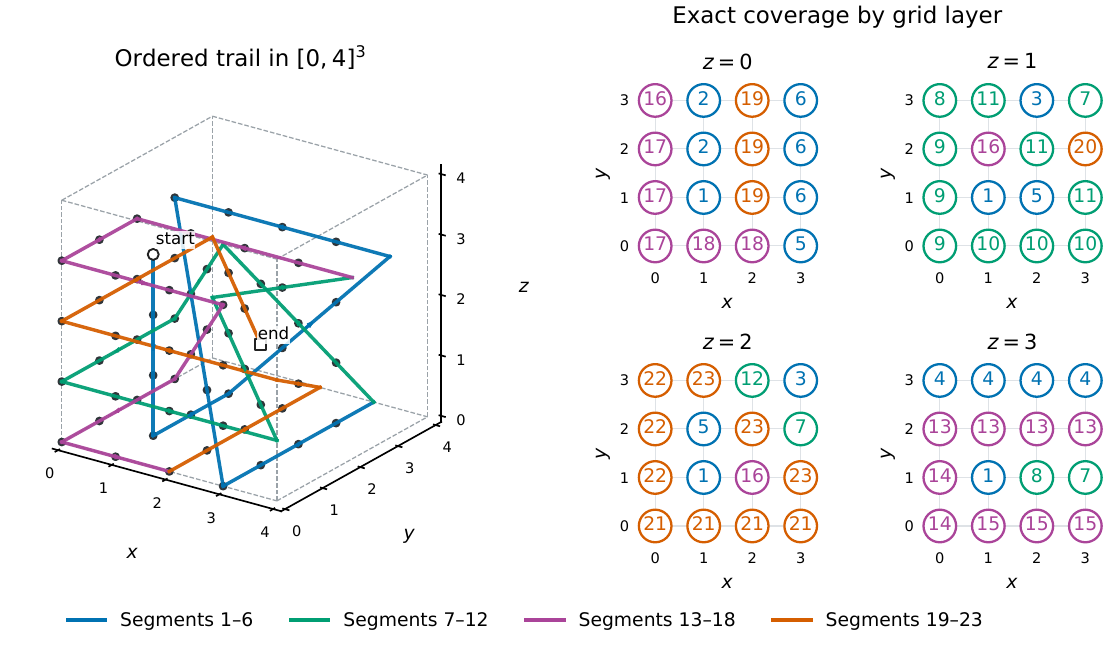}
\captionsetup{singlelinecheck=false}
\caption[An optimal covering trail inside {$[0,4]^3$}]{A $23$-segment trail covering all $64$ points of $P(4,3)$ inside $[0,4]^3$; its ordered vertices are:
\[
\begin{aligned}
&(1,1,3)\to(1,1,0)\to(1,3,0)\to(4,3,3)\to(0,3,3)\to(3,0,0)\\
 &\to(3,4,0)\to(3,0,4)\to(0,3,1)\to(0,0,1)\to(4,0,1)\to(0,4,1)\\
 &\to(4,2,3)\to(0,2,3)\to(0,0,3)\to(3,0,3)\to(0,3,0)\to(0,0,0)\\
 &\to(2,0,0)\to(2,4,0)\to(4,0,2)\to(0,0,2)\to(0,4,2)\to(3,1,2).
\end{aligned}
\]
Left: the full trail, with the grid points in black and the allowed box dashed. Right: the four horizontal grid layers. Each number identifies the first segment containing that point, counted from $1$ in the displayed sequence; colours agree across panels (the figure uses exactly the coordinates certified in Lean and published in \cite{technical}).}\label{fig:boxed23}
\end{figure}

\begin{proof}
Theorem~\ref{thm:main23} gives the lower bound for arbitrary real vertices without any box constraint. To attain it, join the $24$ vertices displayed in Figure~\ref{fig:boxed23} in their stated order. The resulting $23$ segments are nondegenerate, and for each of the $64$ grid points the certificate in \cite{technical} specifies a segment and a rational parameter in $[0,1]$ expressing the point as a convex combination of that segment's endpoints. Exact verification of these incidences proves coverage. The supporting lines are pairwise distinct, which excludes a common subsegment of positive length and hence retracing. Finally, all vertices lie in $[0,4]^3$, so convexity keeps every segment inside the box.
\end{proof}

The integer coordinates of this witness make it easy to describe, but integrality is not an assumption of the constrained problem or its lower bound. Its optimality therefore holds among all trails with real vertices in $[0,4]^3$. Since the construction does use points outside $[0,3]^3$, it does not settle the optimum in that smaller cube.

\subsection{Formal certificates and trust assumptions}\label{sec:certificates}
The geometric reductions in Section~\ref{sec:four} and Subsection~\ref{sec:five} reduce the covering problems for polygonal chains with arbitrary real vertices to finite predicates. The role of the computation is to discharge those predicates, not to replace the reductions by a search over a chosen lattice of vertices. The chain form of the capacity theorem is formalized as \path{Ripa.grid_chain_capacity_all}, for arbitrary finite dimension and $n\ge3$. Its proof is checked by the kernel, including the finite scalar case distinctions in Lemma~\ref{geo:lem:ternary}.

The two numerical results are assembled in \path{Ripa.certified_optimality_four_three_23} and \path{Ripa.certified_lower_bound_five_three_32}, with conclusions
\[
\begin{aligned}
 &\texttt{TrailCoverable 4 3 23}\;\land\;
   \forall r,\;\texttt{Coverable 4 3 r}\to23\le r,\\
 &\forall r,\;\texttt{Coverable 5 3 r}\to32\le r,
\end{aligned}
\]
respectively. Here, \texttt{Coverable} refers to arbitrary polygonal chains, while \texttt{TrailCoverable} imposes Definition~\ref{def:covering}, so these conclusions include the extension in Remark~\ref{rem:walks}. All geometric, graph and packing hypotheses in the final assembly have been discharged. The dependency audits list the standard axioms \texttt{Classical.choice}, \texttt{Quot.sound} and \texttt{propext}, together with native-evaluation dependencies and no \texttt{sorryAx}.

Native reflection additionally trusts the native compiler and runtime \cite{native}; accordingly, the two numerical exclusions are not kernel-only proofs. The recorded runs also reuse compiled library artifacts, and the release includes the transitive source dependencies needed to rebuild them. The technical note \cite{technical} documents the pinned Lean version, exact predicates, file hashes, axiom audits and reproduction procedure. Its supplementary artifact for $P(5,3)$ contains $254$ native-evaluation dependencies.

For the construction in Figure~\ref{fig:boxed23}, the theorem \path{Ripa.box23_covering_trail} checks coverage, nondegeneracy, absence of retracing and box containment in the kernel. The theorem \path{Ripa.box23_optimality_of_unrestricted_lower} then transfers the unrestricted lower bound to the constrained problem; that lower bound retains the native-reflection assumptions just described. The supplementary theorem \path{LargeGridReview.trailCoverable_2000_lower} formalizes the complete counting argument in Appendix~\ref{app:largegrid}, from the original real-vertex definition \path{Ripa.TrailCoverable}. It is checked by the kernel without native reflection; its axiom audit contains only \path{propext}, \path{Quot.sound} and \path{Classical.choice}. The earlier constructions cited in Table~\ref{tab:bounds}, including the general upper bounds for $n=2$ and $n=3$, are results from the literature and are not newly formalized here.

\subsection{Relation to the proposed general formulae}\label{sec:conjecture}
The alternative general formulae in \cite{conjecture} give $21$, $22$ and $23$ for $h(4,3)$. Theorem~\ref{thm:main23} excludes the first two at this grid. For $n\ge3$ and $k\ge2$, the proposed formula agreeing with $h(4,3)=23$ is
\begin{equation}\label{eq:conjecture}
 h(n,k)\stackrel{?}{=}\frac{n^k-1}{n-1}+(k-1)(n-3),
\end{equation}
and it also reproduces the known planar and ternary optima. However, this agreement does not imply that the formula holds throughout its proposed range. Motivated by the heavy-segment counting method of Dumitrescu and T\'oth \cite[Section~2]{dumitrescutoth}, we derive a separate lower bound for polygonal trails without relying on the tree structure considered in their original argument. Appendix~\ref{app:largegrid} establishes this bound directly for trails allowing isolated self-intersections and proves
\[
 h(2000,3)\ge4010357>4005995,
\]
where the right-hand number is the prediction of~\eqref{eq:conjecture}. Thus, none of the three proposed formulae holds throughout its stated range.

For $P(5,3)$, the value $35$ suggested by~\eqref{eq:conjecture} remains compatible with the bounds
\[
 32\le h(5,3)\le36
\]
from Proposition~\ref{prop:five32} and the construction in \cite{small}, but is not established by these results. The exact value for this grid remains open. Our box-constrained construction attains the unrestricted optimum for $P(4,3)$ while staying entirely inside the enlarged box~$[0,4]^3$.

\section*{Declaration of AI assistance}
OpenAI ChatGPT/Codex tools were used extensively by Aldo Roberto Pessolano during the development of this work. Their substantive use included the proof strategies and arguments for the capacity bound in Section~\ref{sec:capacity}, the geometric reductions and finite exclusions in Sections~\ref{sec:four} and~\ref{sec:formal}, programming and computational experiments, the Lean 4 formalization, bibliographic checking, the supplementary arguments in Appendices~\ref{app:ternarycases} and~\ref{app:largegrid}, and preparation and revision of the text. Some mathematical arguments were therefore developed interactively with these systems rather than originating solely from the human authors.

Marco Rip\`a contributed previously published results and constructions forming part of the mathematical basis of this work and participated in the review, formulation, organization and exposition of the present manuscript. The human authors are responsible for checking the mathematical statements and citations and for approving the version submitted for publication.

\appendix

\section{The six exceptional ternary configurations}\label{app:ternarycases}
This appendix supplies the elementary check used in Lemma~\ref{geo:lem:ternary}.

Fix the coordinate indices $\alpha,\beta$ from Lemma~\ref{geo:lem:ternary}, so that $D=d_\alpha e_\beta-d_\beta e_\alpha$ is nonzero and divisible by $3$. The four entries of the corresponding $2\times2$ matrix lie in $\{-2,-1,0,1,2\}$. A zero entry would make the determinant a product of two entries, which cannot be a nonzero multiple of $3$. Thus, all entries are nonzero, and the possible nonzero determinants divisible by $3$ are $D\in\{-6,-3,3,6\}$.

Use the grid points $a,b\in P(3,k)$, directions $d,e$, and real parameters $\lambda,\mu\ge0$ from that lemma. The relevant constraints in the selected coordinates are
\[
\begin{cases}
a_\gamma-\lambda d_\gamma=b_\gamma-\mu e_\gamma,\\
0\le a_\gamma,\ a_\gamma+d_\gamma,\ b_\gamma,\ b_\gamma+e_\gamma\le2
\end{cases}
\quad(\gamma\in\{\alpha,\beta\}).
\]
Reflecting either coordinate ($x_\gamma\mapsto2-x_\gamma$, $\gamma\in\{\alpha,\beta\}$), exchanging the two coordinates, and exchanging the two segments preserve their feasibility.
The reflections change the signs of the rows of the matrix with columns $(d_\alpha,d_\beta)$ and $(e_\alpha,e_\beta)$, while exchanging the coordinates or the two segments exchanges its rows or columns, respectively. 

To count the admissible direction pairs, note that exactly one, two, or three of the four entries have absolute value $2$, and that in the two-entry case these occupy opposite corners. There are $4$, $2$, and $4$ possible placements, respectively. Of the $16$ sign assignments for each placement, exactly eight give $|D|\in\{3,6\}$: the two determinant products must have opposite signs with one or three entries of absolute value $2$, and the same sign with two. Hence, there are $(4+2+4)\cdot8=80$ admissible direction pairs. 

Applying the stated row-sign changes and row and column exchanges to the six representatives in Table~\ref{tab:ternarycases} gives pairwise disjoint orbits of sizes $16,16,8,16,16,8$, whose sum is $80$; they therefore exhaust all admissible pairs. The four orbits with $|D|=3$ contain $48$ pairs, and the two with $|D|=6$ contain $32$. Swapping the two coordinates reverses the sign of $D$, so there are $24$ direction pairs for each of $D=3,-3$ and $16$ for each of $D=6,-6$.

\begin{table}[H]
\centering\small
\renewcommand{\arraystretch}{1.4}
\setlength{\tabcolsep}{5pt}
\begin{tabular}{|c|c|c|c|p{9cm}|}
\hline
$(d_\alpha,d_\beta)$ & $(e_\alpha,e_\beta)$ & $D$ & Orbit & \multicolumn{1}{c|}{\parbox{9cm}{\centering Consequence of the coordinate bounds}}\\
\hline
$(1,1)$ & $(-2,1)$ & $3$ & $16$ & $a_\alpha\le1$, $b_\alpha=2$, so $a_\alpha-\lambda=2+2\mu$ is impossible.\\
$(1,1)$ & $(-1,2)$ & $3$ & $16$ & $a_\alpha\le1\le b_\alpha$, so $\lambda+\mu=a_\alpha-b_\alpha\le0$. Hence $\lambda=\mu=0$.\\
$(1,2)$ & $(-2,-1)$ & $3$ & $8$ & $a_\alpha\le1$, $b_\alpha=2$: the same contradiction as in the first row.\\
$(1,2)$ & $(-2,2)$ & $6$ & $16$ & $a_\alpha\le1$, $b_\alpha=2$: again $a_\alpha-\lambda=2+2\mu$ is impossible.\\
$(1,2)$ & $(2,-2)$ & $-6$ & $16$ & $a_\beta=0$, $b_\beta=2$, so $-2\lambda=2+2\mu$ is impossible.\\
$(1,2)$ & $(2,1)$ & $-3$ & $8$ & $a_\beta=b_\alpha=0$. The two equations give $\lambda-2\mu=a_\alpha$ and $\mu-2\lambda=b_\beta$, so $-\lambda-\mu=a_\alpha+b_\beta\ge0$. Hence $\lambda=\mu=0$.\\
\hline
\end{tabular}
\caption{Exceptional direction pairs in two coordinates, up to signed coordinate permutations and exchange of the segments. Parameters remain real throughout.}\label{tab:ternarycases}
\end{table}
\FloatBarrier
In the two feasible classes, $\lambda=\mu=0$ forces $a=b=v$ in every coordinate of the original configuration. This proves the exceptional implication without using numerical approximation or assuming that the parameters are integers.

\section{A large-grid obstruction for the proposed formula}\label{app:largegrid}

The argument below adapts the heavy-segment and interior-line counting method of Dumitrescu and T\'oth \cite[Section~2]{dumitrescutoth}. We give an explicit estimate for a single grid, using coarse direction counts. The proof counts the original segments of the trail, without subdividing them at crossings.

\begin{proposition}\label{prop:largegrid}
With the definition of trails used in this paper, including permission for isolated self-intersections and repeated visits to grid points,
\[
 h(2000,3)\ge4010357.
\]
\end{proposition}
\begin{proof}
Put $n=2000$ and consider any covering trail with $r$ segments. We use a boundary width of $40$ grid steps. Call a segment \emph{heavy} in this proof if it contains at least $n-40=1960$ grid points, and let $H$ be their number. A nonconstant coordinate is injective on the points of a line, so each segment contains at most $n$ grid points; a segment with at most one grid point satisfies this bound as well. A nonheavy segment contains at most $n-41$ grid points, but we use the weaker bound $n-40$ to simplify the counting argument. Counting incidences, with repeated visits counted more than once, gives
\begin{equation}\label{eq:heavycount}
 n^3\le nH+(n-40)(r-H)=(n-40)r+40H.
\end{equation}

We first recover the possible directions of a heavy support directly from its grid points. Choose a nonconstant coordinate $\alpha$. Its projection $T\subseteq\{0,\ldots,1999\}$ has at least $1960$ elements and contains two consecutive integers: otherwise, $T$ and $T+1$ would be disjoint subsets of $\{0,\ldots,2000\}$ with at least $3920$ elements in total. Two grid points whose $\alpha$-th coordinates differ by one determine an integer direction $d$ with $d_\alpha=1$. Collinearity then gives
\[
 b_\beta=a_\beta+(b_\alpha-a_\alpha)d_\beta
\]
for every grid point $b$ of the segment, where $a$ is the first chosen grid point. The span of $T$ is at least $1959$. If $|d_\beta|\ge2$, the span in coordinate $\beta$ would exceed $1999$, a contradiction. Thus, every component of $d$ lies in $\{-1,0,1\}$.

Distinct heavy segments have distinct supports. Indeed, the union of their grid-point sets would contain at most $n$ points on a common support, and their intersection at most one point by the prohibition on retracing. Their cardinalities would therefore sum to at most $n+1=2001$, whereas two heavy segments contribute at least $3920$.

There are six unoriented directions with exactly two nonzero coordinates. For each, the constant coordinate has $n$ possibilities. After reflecting a varying coordinate if needed, the varying coordinates differ by a constant offset. The projection on either varying coordinate has minimum at most $40$ and maximum at least $1959$; requiring both coordinates to lie between $0$ and $1999$ bounds the offset between $-40$ and $40$. Hence, there are at most $81$ possibilities. These supports therefore number at most $6\cdot81n$.

There are four unoriented directions with three nonzero coordinates. After coordinate reflections, each has direction $(1,1,1)$ and is described by two offsets. The same argument bounds each offset between $-40$ and $40$. Thus, the number of non-axis-aligned heavy supports is at most
\[
 N_{\mathrm{diag}}:=6\cdot81n+4\cdot81^2=998244.
\]

For each axis direction, call a support \emph{interior} when both transverse coordinates belong to $\{40,\ldots,n-41\}$. The number of other supports in that direction is
\[
 N_{\mathrm{bdry}}:=n^2-(n-80)^2=160(n-40)=313600.
\]
This count can also be obtained by splitting the boundary into four disjoint strips, of sizes $40n$, $40n$, $40(n-80)$ and $40(n-80)$. Across the three axis directions, at most $3N_{\mathrm{bdry}}$ heavy supports are not interior. The remaining axis-aligned heavy segments split into three families. Choose a largest family, call its direction vertical, and let $K$ be its cardinality. Then
\begin{equation}\label{eq:interiorheavy}
 3K\ge H-N_{\mathrm{diag}}-3N_{\mathrm{bdry}}.
\end{equation}

Each of these $K$ segments, except possibly the last segment of the trail, has a successor. Let $v$ be its common endpoint with the successor. Both transverse coordinates of $v$ equal the integer coordinates of the heavy support. This follows, for example, by determining the support from two distinct grid points on the heavy segment.

If the successor is nonvertical, choose a transverse coordinate $\beta$ that varies strictly along it. Since $v_\beta\in\{40,\ldots,n-41\}$, its grid points have distinct integer coordinate values on one side of $v_\beta$, including the endpoint value when it occurs. There are at most $n-v_\beta\le n-40$ such values in the increasing direction, or $v_\beta+1\le n-40$ in the decreasing direction. If the successor is vertical, it lies on the same support as the heavy segment. Their grid-point sets have union of size at most $n$ and intersection of size at most one. The successor therefore contains at most $41\le n-40$ grid points. This second case avoids any need to assume that the trail is minimal or to merge collinear segments.

Distinct selected segments have distinct successor positions. Removing the last position, when present, therefore leaves at least $K-1$ distinct segments that each lose at least $40$ points from the individual capacity $n$. Subtracting these deficits from the initial incidence bound $rn$ gives
\begin{equation}\label{eq:successordeficit}
 n^3\le rn-40(K-1).
\end{equation}
For $K=0$ the inequality holds as well. No step uses the absence of crossings between nonconsecutive segments.

Combining~\eqref{eq:interiorheavy} and~\eqref{eq:successordeficit} gives
\[
3n^3\le3rn-40H+40(N_{\mathrm{diag}}+3N_{\mathrm{bdry}}+3).
\]
On the other hand, \eqref{eq:heavycount} implies
\[
40H\geq n^3-(n-40)r.
\]
Substituting this lower bound for $40H$ into the preceding inequality yields
\[
4n^3\le(4n-40)r+40(N_{\mathrm{diag}}+3N_{\mathrm{bdry}}+3),
\]
and therefore
\[
r(4n-40)\ge4n^3-40(N_{\mathrm{diag}}+3N_{\mathrm{bdry}}+3)
=31922438120.
\]
Since $4n-40=7960$ and $r$ is an integer,
\[
 r\ge\left\lceil\frac{31922438120}{7960}\right\rceil
 =4010357.
\]
\end{proof}
For $k=3$, the expression proposed in~\eqref{eq:conjecture} simplifies to
\[
\frac{n^3-1}{n-1}+2(n-3)=n^2+3n-5.
\]
At $n=2000$, it predicts $4005995$ segments, whereas Proposition~\ref{prop:largegrid} requires at least $4010357$. Since the latter is strictly greater than the former, the proposed formula cannot hold for all $n\ge3$ and $k\ge2$.

The proof of Proposition~\ref{prop:largegrid} is formalized by the theorem \path{LargeGridReview.trailCoverable_2000_lower}, whose only mathematical premise is \texttt{Ripa.TrailCoverable 2000 3 r}. This proof uses kernel checking and no native reflection; its axiom audit contains only \texttt{propext}, \texttt{Quot.sound} and \texttt{Classical.choice}.

\end{document}